\documentclass[12pt]{elsart}
\usepackage{amsfonts,amsmath,amssymb}

\newtheorem{df}[thm]{\bf Definition}
\newtheorem{ex}[thm]{\bf Example}

\newcommand{\BZ}{B^{\mathbb Z}}

\newcommand{\TR}{T^{R,b}}

\newcommand{\calC}{\mathcal C}

\journal{}

\date{}

\begin{document}

\begin{frontmatter}
\title
{The type ${\rm III_1}$ factor generated by regular representations
of the infinite dimensional nilpotent group $B_0^\mathbb Z$ }

\thanks[dec]{The authors would like to thank the
Max-Planck-Institute of Mathematics, Bonn  for the hospitality. The
first author was supported by the Max-Planck Doctoral scholarship.
The second author was partially supported by the DFG  project 436
UKR 113/87.}
\author{Ivan Dynov}
\address{
Department of Mathematics and Statistics, York University N536 Ross
Building, \\
4700 Keele Street, Toronto, ON M3J 1P3 Canada\\
Max-Planck-Institut f\"ur Mathematik, Vivatsgasse 7, D-53111 Bonn,
Germany,\\
 E-mail: dynov@mpim-bonn.mpg.de} and
\author{Alexandre Kosyak
\corauthref{cor}} \ead{kosyak01@yahoo.com, kosyak@imath.kiev.ua}
\address{Max-Planck-Institut f\"ur Mathematik, Vivatsgasse 7, D-53111 Bonn,
Germany}
\address{Institute of Mathematics, Ukrainian National Academy of Sciences,
3 Tereshchenkivs'ka, Kyiv, 01601, Ukraine,\\
E-mail: kosyak01@yahoo.com, kosyak@imath.kiev.ua \\
tel.: 38044 2346153 (office), 38044 5656758 (home), fax: 38044
2352010}

\corauth[cor]{Corresponding author}

\newpage

\begin{abstract}
 We study the von Neumann algebra, generated by the regular
representations of the infinite-dimensional nilpotent group
$B_0^{\mathbb Z}$. In \cite{Kos00f} a condition have been found on
the measure for the right von Neumann algebra to be the commutant of
the left one. In the present article, we prove that, in this case,
the von Neumann algebra generated by the regular representations of
group $B_0^{\mathbb Z}$ is the type ${\rm III}_1$ hyperfinite
factor.

We use a technique, developed in \cite{Kos08arIII.1} where a similar
result was proved for the group $B_0^{\mathbb N}$. The crossed
product allows us to remove some technical condition on the measure
used in \cite{Kos08arIII.1}.
\end{abstract}

\begin{keyword}  von Neumann algebra, modular operator, operator of canonical
conjugation, type ${\rm III}_1$ factor, unitary representation,
infinite-dimensional groups, nilpotent groups, regular
representations, irreducibility, infinite tensor products, Gaussian
measures, Ismagilov conjecture

\MSC 22E65 \sep (28C20, 43A80, 58D20)
\end{keyword}

\end{frontmatter}

\sloppy

\newcommand{\tr}{\mathrm{tr}\,}
\newcommand{\rank}{\mathrm{rank}\,}
\newcommand{\diag}[1]{\mathrm{diag}\,(#1)}
\renewcommand{\Im}{\mathrm{Im}\,}

\maketitle
\newpage
\tableofcontents

\section{ Introduction}
 We study the von Neumann algebra, generated by the regular
representations of the infinite-dimensional nilpotent group
$B_0^{\mathbb Z}$. The conditions of irreducibility of the regular
and quasiregular representations of infinite-dimensional groups
(associated with some quasi-invariant Gaussian measures $\mu_b$) are
given by the so-called Ismagilov conjecture (see
\cite{Kos90,Kos92,Kos94}). In this case the corresponding von
Neumann algebra is a type ${\rm I}_\infty$ factor. In \cite{Kos00f}
a condition have been found on the measure for the right von Neumann
algebra to be the commutant of the left one.

In the present article, we prove that, in this case, the von Neumann
algebra generated by the regular representations of the
infinite-dimensional nilpotent group $B_0^{\mathbb Z}$ is the type
${\rm III}_1$ hyperfinite factor. Moreover this factor is unique.
\newpage
We recall that the {\it first examples of a non-type I factor},
namely a type ${\rm II}_1$ factor, were also obtained by Murray and
von Neumann as von Neumann algebra {\it generated by the regular
representation} of a discrete ICC group (i.e. the group for which
all conjugacy classes are infinite, except the trivial one)
We shall show that the regular representations of non-discrete
infinite-dimensional  groups provide examples of non-type I or II,
but type III factors, namely the type ${\rm III}_1$.
\section{ Regular representations}
Let us consider the group $\tilde G=B^{\mathbb Z}$ of all
upper-triangular real matrices of infinite order with units on the
diagonal
$$
\tilde G=B^{\mathbb Z}=\{I+x\mid x=\sum_{k,n\in{\mathbb Z}\,
k<n}x_{kn}E_{kn}\},
$$
and its subgroup
$$
G=B_0^{\mathbb Z}=\{I+x\in B^{\mathbb Z}\mid\,\,x\,\,{\rm
is}\,\,{\rm finite}\},
$$
where $E_{kn}$ is an infinite-dimensional matrix with $1$ at the
place $k,n\in {\mathbb Z}$ and zeros elsewhere,
 $x=(x_{kn})_{k<n}$ is {\it finite} means that $x_{kn}=0$ for
all $(k,n)$ except for a finite number of indices $k,n\in{\mathbb
Z}$.

Obviously, $B_0^{\mathbb Z}=\varinjlim_{n} B(2n-1,{\mathbb R})$ is
the inductive limit of the group $B(2n-1,{\mathbb R})$ of real
upper-triangular matrices with units on the principal diagonal
realized in the following form
$$
B(2n-1,{\mathbb R})=\{I+\sum_{-n+1\leq k<r\leq n-1 }x_{kr}E_{kr}\mid
x_{kr}\in {\mathbb R}\},\quad n\in {\mathbb N},
$$
with respect to the symmetric imbedding
$$B(2n-1,{\mathbb R})\ni
x\mapsto i^s(x)=x+E_{-n,-n}+E_{nn}\in B(2n+1,{\mathbb R}).$$

We define the Gaussian measure $\mu_b$ on the group $B^{\mathbb Z}$
in the following way
\begin{equation}
\label{mu-bZ} d\mu_b(x) = \otimes_{ k,n\in{\mathbb Z},\,
k<n}(b_{kn}/\pi)^{1/2}\exp(-b_{kn}x_{kn}^2 )dx_{kn}=
\otimes_{k,n\in{\mathbb Z},\,k<n}d\mu_{b_{kn}}(x_{kn}),
\end{equation}
where $b=(b_{kn})_{k<n}$ is some set of positive numbers
$b_{kn}>0,\,\,k,n\in{\mathbb Z}$.

Let us denote by $R$ and $L$ the right and the left action of the
group $ B^{\mathbb Z}$ on itself: $R_t(s)=st^{-1},\,\,\,L_t(s)= ts,
\,\,s,t \in B^{\mathbb Z}$ and by $ \Phi:B^{\mathbb Z}\mapsto
B^{\mathbb Z},\,\,\Phi(I+x):=(I+x)^{-1}$ the inverse mapping.  It is
known \cite{Kos01m2,Kos04} that
\begin{lem}
$\mu_b^{R_t}\sim \mu_b\,\,\forall t\in B_0^{\mathbb Z}$ if and only
if $S^R_{kn}(b)<\infty,\,\,\forall k,n\in{\mathbb Z},\,\,k<n$ where
\begin{equation} \label{SRkn}
S^R_{kn}(b)=\sum_{r=-\infty}^{k-1}\frac{b_{rn}}{b_{rk}}.
\end{equation}
\end{lem}
\begin{lem}
$\mu_b^{L_t}\sim \mu_b\,\,\forall t\in B_0^{\mathbb Z}$ if and only
if $S^L_{kn}(b)<\infty,\,\,\forall k,n\in{\mathbb Z},\,\, k<n$,
where
\begin{equation}
 \label{SLkn}
S^{L}_{kn}(b)=\sum_{m=n+1}^\infty\frac{b_{km}}{b_{nm}}.
\end{equation}
\end{lem}

\begin{lem}
$\mu_b^{L_{I+tE_{kn}}}\perp \mu_b\,\, \forall t
\in {\mathbb R}\backslash\{0\} \Leftrightarrow
S^{L}_{kn}(b)=\infty, \,\,\, 
k,n\in{\mathbb Z},\,\, k<n.$
\end{lem}
Let us denote
\begin{equation}
\label{E(b)}
 E(b)=\sum_{k<n<r}\frac{b_{kr}}{b_{kn}b_{nr}},\,\,
E_m(b)=\sum_{k<n<r\leq
m}\frac{b_{kr}}{b_{kn}b_{nr}},\,\,m\in{\mathbb Z}.
\end{equation}
\begin{lem}{\rm \cite{KosZek00}}
If $E(b)
<\infty$, then $\mu_b^\Phi\sim\mu_b.$
\end{lem}

\begin{rem}{\rm \cite{KosZek00}}
If $\mu_b^\Phi\sim\mu_b$ then $\mu_b^{L_t}\sim \mu_b\Leftrightarrow
\mu_b^{R_t}\sim \mu_b\,\,\forall t\in B_0^{\mathbb Z}.$
\end{rem}
\begin{pf}
This follows from the fact that the inversion $\Phi$ replace the
right and the left action: $R_t\circ\Phi=\Phi\circ L_t\,\,\forall t
\in B^{\mathbb Z}$. Indeed, if we denote $\mu^{f}(C)=\mu(f^{-1}(C))$
for a measurable set $C$, we have $(\mu^{f})^g=\mu^{f\circ g}$.
Hence
$$
\mu_b\sim\mu_b^{R_t}\sim(\mu_b^{R_t})^\Phi= \mu_b^{R_t\circ\Phi}=
\mu_b^{\Phi\circ L_t} =(\mu_b^\Phi)^{L_t}\sim\mu_b^{L_t},\,\,\forall
t \in B_0^{\mathbb Z}.
$$
\qed\end{pf}
\begin{rem} We have
\begin{equation}
\label{}
 E(b)
=\sum_{k<n}\frac{S^{L}_{kn}(b)}{b_{kn}}=\sum_{k<n}\frac{S^{R}_{kn}(b)}{b_{kn}}
,\,\, E_m(b)=\sum_{k<n\leq m}\frac{S^{R}_{kn}(b)}{b_{kn}}.
\end{equation}
 Indeed
$$
\sum_{k<n}\frac{S^{L}_{kn}(b)}{b_{kn}}=\sum_{k<n}\sum_{r=n+1}^\infty
\frac{b_{kr}}{b_{kn}b_{nr}}=
\sum_{k<n<r}\frac{b_{kr}}{b_{kn}b_{nr}}=E(b)
$$
$$
=\sum_{n<r}\frac{1}{b_{nr}}\sum_{k=-\infty}^{n-1}\frac{b_{kr}}{b_{kn}}=
\sum_{n<r}\frac{S^{R}_{nr}(b)}{b_{nr}}.
$$
\end{rem}
If $\mu_b^{R_t}\sim \mu_b$ and $\mu_b^{L_t}\sim \mu_b\,\,\forall
t\in B_0^{\mathbb Z}$, one can define in a natural way (see
\cite{Kos90,Kos92}), an analogue of the right $T^{R,b}$ and the left
$T^{L,b}$ regular representations of the group $B_0^{\mathbb Z}$ in
the Hilbert space $H_b=L^2(B^{\mathbb Z},\mu_b)$
$$
T^{R,b},\,\,T^{L,b}:B_0^{\mathbb Z}\rightarrow U(H_b=L^2(B^{\mathbb
Z},\mu_b)),
$$
$$(T^{R,b}_t f)(x)=(d\mu_b(xt)/d\mu_b(x))^{1/2}f(xt),$$
$$
(T^{L,b}_s f)(x)=(d\mu_b(s^{-1}x)/d\mu_b(x))^{1/2}f(s^{-1}x).
$$
\section{ Von Neumann  algebras}
Let  ${\mathfrak A}^{R,b}=(T^{R,b}_t \mid t \in B_0^{\mathbb
Z})^{\prime\prime}$ (resp.  ${\mathfrak A}^{L,b}=(T^{L,b}_s \mid s
\in B_0^{\mathbb Z})^{\prime\prime}$) be the von Neumann algebras
generated by the right $T^{R,b}$ (resp. the left $T^{L,b}$) regular
representation of the group $B_0^{\mathbb Z}$. 
\begin{thm}
\label{com-thm}{\rm \cite{Kos00f}} (The commutation theorem) If
$E(b)<\infty$ then $\mu_b^\Phi\sim\mu_b$. In this case the right and
the left regular representations are well defined and the
commutation theorem holds:
\begin{equation}
\label{(R)'=(L)-Z} ({\mathfrak A}^{R,b})'={\mathfrak A}^{L,b}.
\end{equation}
 Moreover, the operator $J_{\mu_b}$ given by
\begin{equation}
\label{J-Z}
(J_{\mu_b}f)(x)=(d\mu_b(x^{-1})/d\mu_b(x))^{1/2}\overline{f(x^{-1})}
\end{equation}
is an intertwining operator:
$$
T^{L,b}_t=J_{\mu_b}T^{R,b}_tJ_{\mu_b},\,\,t\in B_0^{\mathbb
Z}\,\,\text{\, and\,}\,\, J_{\mu_b}{\mathfrak
A}^{R,b}J_{\mu_b}={\mathfrak A}^{L,b}.$$
\end{thm}
If $\mu_b^{R_t}\sim \mu_b\,\,\forall t\in B_0^{\mathbb Z}$ but
 $\mu_b^{L_t}\perp \mu_b\,\,\forall t\in B_0^{\mathbb
Z}\backslash\{e\}$ one can't define the left regular representation
of the group $B_0^{\mathbb Z}$. Moreover the following theorem holds
\begin{thm}
\label{irr.07II} {\rm \cite{Kos01m2}} The right regular
representation $T^{R,b}:B_0^{\mathbb Z}\mapsto
U(H_b)$ is irreducible if\\
(1) $\mu_{b}^{L_{s}}\perp
\mu_{b}\,\, \forall s\in B_0^{\mathbb Z}\backslash\{0\},$\\
(2) the measure $\mu_b$ is $B_0^{\mathbb Z}$ right-ergodic,\\
(3) $\sigma_{kn}(b)=\infty,\,\, \forall  k<n,\,k,n\in{\mathbb Z}$,
where
$$
\sigma_{kn}(b)=\sum_{m=n+1}^\infty
\frac{b^2_{km}}{[S^{R}_{km}(b)+b_{km}][S^{R}_{nm}(b)+b_{nm}]}.
$$
\end{thm}
\begin{rem}
We do not know whether the Ismagilov conjecture holds in this case,
namely, whether conditions 1) and 2) of the theorem are the criteria
of the irreducibility of the representation $T^{R,b}$ of the group
$B_0^{\mathbb Z}$ as holds for example for the group $B_0^{\mathbb
N}$ (see \cite{Kos90,Kos92}).
\end{rem}
\begin{rem}\label{rem_ergo}
We do not know the criterion of the $B_0^{\mathbb Z}$-ergodicity of
the measure $\mu_b$ on the space $B^{\mathbb Z}$. The sufficient
conditions are the following $E_m(b)<\infty$  for all $m\in{\mathbb
Z}$.
\end{rem}
\begin{rem}
The von Neumann algebra ${\mathfrak A}^{R,b}$ is a type $I_\infty$
factor if the conditions of the Theorem \ref{irr.07II} are valid.
\end{rem}
Let us  assume now that $\mu_b^{L_t}\sim
\mu_b\sim\mu_b^{R_t}\,\,\forall t\in B_0^{\mathbb Z}$. In this case
the right regular representation and the left regular representation
of the group $B_0^{\mathbb Z}$ are well defined.

In \cite{KosZek00} the condition were studied  {\it when the von
Neumann algebra ${\mathfrak A}^{R,b}$ is a factor}, i.e.
$${\mathfrak A}^{R,b}\cap({\mathfrak A}^{R,b})^\prime= \{\lambda
 {\bf I}\vert\lambda\in {\mathbb C}\}.
$$
Since $T^{L,b}_t \in ({\mathfrak A}^{R,b})^\prime \,\, \forall t\in
B_0^{\mathbb Z}$, we have  ${\mathfrak A}^{L,b}\subset({\mathfrak
A}^{R,b})^\prime$, hence
\begin{equation}
{\mathfrak A}^{R,b}\cap({\mathfrak A}^{R,b})^\prime \subset
({\mathfrak A}^{L,b})^\prime\cap({\mathfrak A}^{R,b})^\prime =
({\mathfrak A}^{R,b}\cup{\mathfrak A}^{L,b})^\prime .
\end{equation}
The last relation shows that ${\mathfrak A}^{R,b}$ is factor if
 the representation
$$
B_0^{\mathbb Z}\times B_0^{\mathbb Z}\ni (t,\,s) \rightarrow
T^{R,b}_t T^{L,b}_s\in U(H_b)
$$
is irreducible. Let us denote
\begin{equation}
\label{S^(R,L)_(kn)(b)} S^{R,L}_{kn}(b)=\sum_{m=n+1}^\infty
\frac{b^2_{km}}{[S^{R}_{km}(b)+b_{km}][S^{L}_{nm}(b)+S^{R}_{nm}(b)]},\,\,\,
k<n.
\end{equation}
\begin{thm}{\rm \cite{KosZek00}} The representation
$$
B_0^{\mathbb Z}\times B_0^{\mathbb Z}\ni (t,\,s) \rightarrow
T^{R,b}_t T^{L,b}_s\in U(H_b)
$$
is irreducible if $S^{R,L}_{kn}(b)=\infty,\,\,\forall k<n$ and the
measure $\mu_b$ is $B_0^{\mathbb Z}$ right-ergodic.
\end{thm}
\begin{cor}\label{cor_factor}
The von Neumann algebra ${\mathfrak A}^{R,b}$ is factor if
$S^{R,L}_{kn}(b)=\infty\,\forall k<n$ and the measure $\mu_b$ is
$B_0^{\mathbb Z}$ right-ergodic.
\end{cor}
\begin{rem}\label{E-inf-B^Z-fact}
In what follows, we shall show that the condition $E(b)<\infty$ is
already sufficient for ${\mathfrak A}^{R,b}$ (and ${\mathfrak
A}^{L,b}$) to be a factor.
\end{rem}
\section{Examples}
In this section we give an example of a measure
$\mu_b,\,\,b=(b_{kn})_{k<n}$ for which $E(b)<\infty$, hence  the
representations $T^{R,b}$ and $T^{L,b}$ are well defined and the
commutation theorem (Theorem \ref{com-thm}) for von Neumann algebras
${\mathfrak A}^{R,b}$ and ${\mathfrak A}^{L,b}$ holds. We show that
the set $b=(b_{kn})_{k<n}$ for which
\begin{equation}
\label{factor?-Z} E(b)<\infty\,\, \text{ and hence }
S^{R}_{kn}(b)<\infty,\,\,S^{L}_{kn}(b)<\infty
\end{equation}
defined respectively by (\ref{E(b)}), (\ref{SRkn}) and (\ref{SLkn}),
is not empty.

In the example (\ref{ex1}) below for the particular case
$b_{kn}=(a_k)^n$ we give some sufficient conditions on the sequence
$a_n$ implying conditions (\ref{factor?-Z}).
\begin{ex}
\label{ex1}
 Let us take $b_{kn}=(a_k)^n,\,\,k,n\in{\mathbb
Z}$.
\end{ex}
We have
\begin{equation}
\label{s^r}
S^R_{kn}(b)=\sum_{r=-\infty}^{k-1}\frac{b_{rn}}{b_{rk}}=\sum_{r=-\infty}^{k-1}
a_r^{n-k}<\infty \,\,\text{\,if\,}\,\sum_{r=-\infty}^0a_r<\infty,
\end{equation}
\begin{equation}
\label{s^l} S^{L}_{kn}(b)=\sum_{m=n+1}^\infty
\left(\frac{a_k}{a_n}\right)^m=
\left(\frac{a_k}{a_n}\right)^{n+1}\sum_{m=0}^\infty
\left(\frac{a_k}{a_n}\right)^m=
\left(\frac{a_k}{a_n}\right)^{n+1}\frac{1}{1-\frac{a_k}{a_n}}<\infty,
\end{equation}
iff $a_k<a_{k+1},\,\,k\in {\mathbb Z}$. Finally we get
$$
E(b)=\sum_{k=-\infty}^\infty\sum_{n=k+1}^\infty\frac{S^{L}_{kn}(b)}{b_{kn}}=
\sum_{k=-\infty}^\infty\sum_{n=k+1}^\infty
\left(\frac{a_k}{a_n}\right)^{n+1}\frac{1}{1-\frac{a_k}{a_n}}\frac{1}{a_k^n}=
$$
$$
\sum_{k=-\infty}^\infty a_k \sum_{n=k+1}^\infty
\left(\frac{1}{a_n}\right)^{n+1}\frac{1}{1-\frac{a_k}{a_n}}
<\sum_{k=-\infty}^\infty
\frac{a_k}{1-\frac{a_k}{a_{k+1}}}\sum_{n=k+1}^\infty
\left(\frac{1}{a_n}\right)^{n+1}
$$
$$
<\sum_{k=-\infty}^\infty
\frac{a_k}{1-\frac{a_k}{a_{k+1}}}\sum_{n=k+1}^\infty
\left(\frac{1}{a_{k+1}}\right)^{n+1}
=\sum_{k=-\infty}^\infty\frac{a_k}{1-\frac{a_k}{a_{k+1}}}
\left(\frac{1}{a_{k+1}}\right)^{k+2} \frac{1}{1-\frac{1}{a_{k+1}}}=
$$
$$
\sum_{k=-\infty}^\infty\frac{\frac{a_k}{a_{k+1}}}
{\left(1-\frac{a_k}{a_{k+1}}\right)^2}
\left(\frac{1}{a_{k+1}}\right)^{k+1}.
$$
\begin{ex}
\label{ex2}
Let us take $b_{kn}=(a_k)^n,\,\,k,n\in{\mathbb Z}$ where
$a_k=s^k,\,\,k\in {\mathbb Z}$ with $s>1$.
\end{ex}
Conditions  (\ref{factor?-Z}) hold for $a_k=s^k$. By (\ref{s^r}) and
(\ref{s^l})  we have
$$
S^{L}_{kn}(b)=\left(\frac{a_k}{a_n}\right)^{n+1}\frac{1}{1-\frac{a_k}{a_n}}=
\left(\frac{1}{s^{n-k}}\right)^{n+1}\frac{1}{1-\frac{1}{s^{n-k}}}\sim
\left(\frac{1}{s^{n-k}}\right)^{n+1},
$$
$$
S^R_{kn}(b)=\sum_{r=-\infty}^{k-1}
a_r^{n-k}=\sum_{r=-\infty}^{k-1}s^{r(n-k)}=
\sum_{r=1-k}^{\infty}\frac{1}{s^{r(n-k)}}=
$$
$$
\left(\frac{1}{s^{n-k}}\right)^{1-k}\frac{1}{1-\frac{1}{s^{n-k}}}\sim
s^{(n-k)(k-1)},
$$
 since
$$
1< \frac{1}{1-\frac{1}{s^{n-k}}}<\frac{1}{1-\frac{1}{s}}.
$$
Using the latter equivalence we conclude that $E(b)<\infty.$ Indeed
we have
$$
E(b)=\sum_{k=-\infty}^\infty\sum_{n=k+1}^\infty\frac{S^{L}_{kn}(b)}{b_{kn}}\sim
\sum_{k=-\infty}^\infty\sum_{n=k+1}^\infty\frac{1}{s^{(n-k)(n+1)}}\frac{1}{s^{kn}}
=\sum_{k=-\infty}^\infty s^k\sum_{n=k+1}^\infty\frac{1}{s^{n(n+1)}}
$$
$$
=\sum_{k=-\infty}^0 s^k\sum_{n=k+1}^\infty\frac{1}{s^{n(n+1)}}+
\sum_{k=1}^\infty s^k\sum_{n=k+1}^\infty\frac{1}{s^{n(n+1)}}<
\sum_{k=-\infty}^0 s^k\sum_{n=-\infty}^\infty\frac{1}{s^{n(n+1)}}
$$
$$
+\sum_{k=1}^\infty
\frac{s^k}{s^{(k+1)^2}}\sum_{n=k+1}^\infty\frac{1}{s^n}
=\sum_{k=-\infty}^0 s^k\sum_{n=-\infty}^\infty\frac{1}{s^{n(n+1)}}+
\sum_{k=1}^\infty \frac{1}{s^{(k+1)^2+1}}<\infty.
$$
\section{Cyclycity}
 We prove that the function $1\in L^2(B^{\mathbb Z},\mu_b)$ is
cyclic and separating for ${\mathfrak A}^{R,b}$ and ${\mathfrak
A}^{L,b}$ if $E(b)<\infty$. We use a method similar to the proof of
ergodicity of $\mu_b$ (under the same condition) in \cite{Kos01m2},
Lemma 4 by reducing the situation to the case of $B^{\mathbb N}$
(\cite{Kos08arIII.1}).
\begin{lem}
If $E(b)<\infty$ then the function $1\in L^2(B^{\mathbb Z},\mu_b)$
is cyclic and separating for $ {\mathfrak A}^{R,b}$.
\end{lem}
\begin{pf}
First we prove that $1$ is cyclic for ${\mathfrak A}^{R,b}$. For any
$m\in\mathbb Z$ we define the subgroups $B^m$ and $B_{(m)}$ of the
group $B^{\mathbb Z}$ as follows:
\begin{equation*}
        B^m:=\{1+x\in B^{\mathbb Z}\mid x=\sum_{k<n\leq m}x_{kn}E_{kn}\},
\end{equation*}
\begin{equation*}
        B_{(m)}:=\{1+x\in B^{\mathbb Z}\mid x=\sum_{k<n,n>m}x_{kn}E_{kn}\}.
\end{equation*}
Obviously, $B^{\mathbb Z}$ is a semi-direct product of the two
groups above ($B_{(m)}$ is a normal subgroup of $B^{\mathbb Z}$) for
any $m$:
\begin{equation*}
B^{\mathbb Z}=B_{(m)}\rtimes B^m,\quad z=yx,\,\,z\in B^{\mathbb
Z},\,\, y\in B_{(m)},\,\,x\in B^m.
\end{equation*}
Let $\mu_{b,(m)}, \,\,\mu_b^m$ be the projections of the measure
$\mu_b$ on the above groups:
\begin{equation*}
\mu_b^m:=\otimes_{k<n\leq m}\mu_{b_{kn}},\quad  \mu_{b,(m)}:
=\otimes_{k<n,n>m}\mu_{b_{kn}}.
\end{equation*}
Then we have
\begin{equation*}
L^2(B^{\mathbb Z},\mu_b)=L^2(B_{(m)},\mu_{b,(m)})\otimes
L^2(B^m,\mu_b^m),\quad f(z)=f(yx).
\end{equation*}
Furthermore let $B_{0,(m)},\,\,B_0^m$ be the intersection of the
above groups with $B_0^{\mathbb Z}$. Now, fix an $m\in\mathbb Z$ and
consider a function $f(z)=f(yx)\in L^2(B^{\mathbb Z},\mu_b)$.
Further, suppose that
\begin{equation}\label{eq_fzero}
(f,a1)=0,\quad \forall a\in{\mathfrak A}^{R,b}.
\end{equation}
First we note that the points of $B_{(m)}$ are invariant under the
right action $R_t$ for all $t\in B_0^m$. Indeed, we have for $t\in
B_0^m$
\begin{equation*}
        (xt)_{kn}=\sum_{j=k+1}^{n-1} x_{kj}t_{jn}=x_{kn}, \quad n>m,
\end{equation*}
since $t_{kn}=\delta_{kn}$ for $n>m$.  We have for $t\in B_0^m$
\begin{align*}
0=&(f,T^{R,b}_t1)=\int_{ B_{(m)}}\int_{B^m} f(yx)T^{R,b}_t1(yx)d\mu_{b,(m)}(y)d\mu^m_b(x)\\
=&\int_{B^m} f_m(x)T^{R,b}_t1(x)d\mu_{b,(m)}(x),
\end{align*}
where
\begin{equation*}
    f_m(x):=\int_{B_{(m)}} f(yx)d\mu_{b,(m)}(y).
\end{equation*}
For the function $f_m$ holds
\begin{equation}\label{eq_fmzero}
    (f_m,T^{R,b}_t1)_m=0,\quad \forall t\in B_0^m,
\end{equation}
where $(.,.)_m$ denotes the restriction of the inner product $(.,.)$
to $L^2(B^m,\mu_b^m)$. Next, we define a bijection $\Psi:B^m \mapsto
B_m$, where $B_m$ is the group
\begin{equation*}
    B_m:=\{1+x; x=\sum_{m\leq k<n}x_{kn}E_{kn}\}\cong B^{\mathbb N},
\end{equation*}
\begin{equation*}
    x^\prime_{kn}=(\Psi(x))_{kn}:=x_{2m-n2m-k}.
\end{equation*}
Note that $\Psi$ are reflections around the axis $k+n=2m$ and if
$m=0$, $x^\prime_{kn}=x_{-n-k}$.  Now we continue with the equation
(\ref{eq_fmzero}):
\begin{eqnarray*}
0=(f_m,T^{R,b}_t1)_m&=&
\int_{B^m}f_m(x)\sqrt{\frac{d\mu^m_b(xt)}{d\mu^m_b(x)}}d\mu_b^m(x)\\
&=&\int_{B_m}f_m^\Psi(x^\prime)\sqrt{\frac{d\mu^{m,\Psi}_b(t^\prime
x^\prime)}{d\mu^{m,\Psi}_b(x^\prime)}}d\mu_b^{m,\Psi}(x^\prime),
\end{eqnarray*}
where $f_m^\Psi:=f_m\circ\Psi$ and
$\mu_b^{m,\Psi}(C)=\mu_b^m(\Psi(C))$ for each Borel set $C$. Since
this holds for all $t\in B_0^m$ and hence all $t^\prime\in B_{0,m}$
and $B_{m}\cong B^{\mathbb N}$,
\begin{equation*}
    0=\int_{B^{\mathbb N}}f^\Psi_m(x)T^{L,b}_t1d\mu_b(x)=(f^\Psi,T^{L,b}_t1)_{\mathbb N},
\end{equation*}
where, more precisely, $f_m^\Psi$ is interpreted as its image under
the isomorphism
 form $B_m$ to $B^{\mathbb N}$ and $(.,.)_{\mathbb N}$ is the inner product on
 $L^2(B^{\mathbb N},\mu_b)$.
It also follows (after taking the linear span and weak limits) that
\begin{equation*}
    (f_m^\Psi,a1)_{\mathbb N}=0,\quad \forall a\in\mathfrak A^{L,b,\mathbb N},
\end{equation*}
where $\mathfrak A^{L,b,\mathbb N}$ are the algebras generated by
the left regular representation $T^{L,b}$ of the group $B_0^{\mathbb
Z}$. But $1$ is cyclic for $\mathfrak A^{L,b,\mathbb N}$, by
\cite{Kos08arIII.1} and hence $f_m^{\Psi}(x^\prime)=0$ for all
$x^\prime\in B_m$. Since $\Psi$ is a bijection, hence we get $f_m=0$
for any $m$.

In addition we have $f_m\rightarrow f$, when $m\rightarrow \infty$
in $L^2(B^{\mathbb Z},\mu_b)$ (see \cite{Kos01m2}, Corollary 1).
 Thus $f_m=0$ for all $m\in\mathbb Z$ implies $f=0$. Since $f$ is arbitrary,
using (\ref{eq_fzero}) we conclude,  that the set ${\mathfrak
A}^{R,b}{\bf 1}$ is dense in $L^2(B^{\mathbb Z},\mu_b)$ and hence
$1$ is cyclic for ${\mathfrak A}^{R,b}$.


Now we turn to the separating property. We know that $1$ is cyclic
for ${\mathfrak A}^{R,b}$. We prove that the same holds for
$({\mathfrak A}^{R,b})^\prime={\mathfrak A}^{L,b}$. Thus, again
consider $f\in L^2(B^{\mathbb Z},\mu_b)$ and assume
\begin{equation}
    (f,b1)=0, \forall b\in{\mathfrak A}^{L,b}.
\end{equation}
Recall that $E(b)<\infty$ implies the existence of the intertwining
operator $J=J_{\mu_b}$, which is anti-unitary. Then the following
calculation holds:
\begin{eqnarray*}
    (f,T^{R,b}_t1)&=&(JT^{R,b}_t1,Jf)\\ \\
              &=&\int\sqrt{\frac{d\mu_b(x^{-1})}{d\mu_b(x)}}
              \sqrt{\frac{d\mu_b((xt)^{-1})}{d\mu_b(x^{-1})}}
              \sqrt{\frac{d\mu_b(x^{-1})}{d\mu_b(x)}}\overline{f(x^{-1})}d\mu_b(x)
              \\ \\
            &=&\int\sqrt{\frac{d\mu_b(t^{-1}x^{-1})}{d\mu_b(x^{-1})}}
            \overline{f(x^{-1})}d\mu(x^{-1}).
\end{eqnarray*}
If we replace $x^{-1}$ by $x$ in the above integral we obtain
$(f,T^{R,b}_t1)=(f,T^{L,b}_t1)$ for all $t\in B_0^{\mathbb Z}$. From
(1) we know that $(f,T^{R,b}_t1)=0$ for all $t\in B_0^{\mathbb Z}$
implies that $f=0$. Hence $(f,T^{L,b}_t1)=0$ for all $t\in
B_0^{\mathbb Z}$ also implies that $f=0$ and hence $1$ is cyclic for
${\mathfrak A}^{L,b}$, since we chose $f$ arbitrarily. \qed\end{pf}
\section{Modular operator}
In this section we recall the construction of the modular operator
for a locally compact group and generalize it to the
infinite-dimensional case. We recall \cite{Dix69C} (see also
\cite{Kos08arIII.1}) how to find the modular operator and the
operator of canonical conjugation for the von Neumann algebra
${\mathfrak A}^\rho_G$, generated by the right regular
representation $\rho$ of a locally compact Lie group $G$. Let $h$ be
a right invariant Haar measure on $G$ and
$$
\rho,\lambda:G\mapsto U(L^2(G,h))
$$
be the right and the left regular representations of the group $G$
defined by
$$
(\rho_tf)(x)=f(xt),\,\,(\lambda_tf)(x)=(dh(t^{-1}x)/dh(x))^{1/2}f(t^{-1}x).
$$
To define  the {\it right Hilbert algebra} on $G$ we can proceed as
follows. Let $M(G)$ be an algebra  of all probability measures on
$G$ with convolution $\mu*\nu$ determined by
\index{algebra!Hilbert!right }
$$
\int_Gf(u)d(\mu*\nu)(u)=\int_G\int_Gf(st)d\mu(s)d\nu(t).
$$
We define the homomorphism
$$
M(G)\ni\mu\mapsto \rho^\mu=\int_G\rho_t d\mu(t)\in B(L^2(G,h)).
$$
We have $\rho^{\mu}\rho^{\nu}=\rho^{\mu*\nu}$, indeed
$$
\rho^{\mu}\rho^{\nu}=\int_G\rho_t d\mu(t)\int_G\rho_s d\nu(s)=
\int_G\int_G\rho_{ts} d\mu(t)\nu(s)=\int_G\rho_t
d(\mu*\nu)(t)=\rho^{\mu*\nu}.
$$
Let us consider a subalgebra $M_h(G):=\{\mu\in M(G)\mid \nu\sim h\}$
of the algebra $M_h(G)$.  In the case when $\mu\in M_h(G)$ we can
associate with the measure $\mu$ its Rodon-Nikodim derivative
$d\mu(t)/dh(t)=f(t)$. When $f\in C_0^\infty(G)$ or $f\in L^1(G)$ we
can write
$$
\rho^f=\int_Gf(t)\rho_tdh(t),
$$
hence we can replace the algebra $M_h(G)$ by its subalgebra
identified with an algebra of functions $C_0^\infty(G)$ or
$L^1(G,h)$ with convolutions.

If we replace the Haar measure $h$ with some  measure $\mu\in
M_h(G)$ we obtain the isomorphic image $T^{R,\mu}$ of the right
regular representation $\rho$ in the space $L^2(G,\mu)$:
$T^{R,\mu}_t=U\rho_tU^{-1}$ where $U:L^2(G,h)\mapsto L^2(G,\mu)$
defined by $(Uf)(x)=\left(\frac{dh(x)}{d\mu(x)}\right)^{1/2}f(x)$.
We have
$$(T^{R,\mu}_tf)(x)=\left(\frac{d\mu(xt)}{d\mu(x)}\right)^{1/2}f(xt),$$
and
$$
T^f=\int_Gf(t)T^{R,\mu}_td\mu(t).
$$
We have (see \cite{{Con94}}, p.462) (we shall write $T_t$ instead of
$T^{R,\mu}_t$ )
$$
S(T^f):=(T^f)^*=\int_G\overline{f(t)}T_{t^{-1}}d\mu(t)=
\int_G\overline{f(t)}T_{t^{-1}}\frac{d\mu(t)}{d\mu(t^{-1})}d\mu(t^{-1})
$$
$$
=\int_G\frac{d\mu(t^{-1})}{d\mu(t)}\overline{f(t^{-1})}T_{t}d\mu(t).
$$
Hence
\begin{equation}
\label{S-n}
(Sf)(t)=\frac{d\mu(t^{-1})}{d\mu(t)}\overline{f(t^{-1})}.
\end{equation}
 In the {\it Tomita-Takesaki theory} \cite{Tak02},
Chapter VI, Lemma 1.2, the operator $S$ is defined for the von
Neumann algebra $M$ of operators on the Hilbert space $H$ by
$$H\ni x\omega\mapsto Sx\omega=x^*\omega\in H,$$ where $x\in M$
and $\omega\in H$ is cyclic (generating) and  separating vector.
\index{theory!Tomita-Takesaki} \index{vector!cyclic}
 \index{vector!separating}

To calculate $S^*$ we use the fact that $S$ is antilinear, i.e.
$(Sf,g)=(S^*g,f)$. We have
$$
(Sf,g)=\int_G\frac{d\mu(t^{-1})}{d\mu(t)}\overline{f(t^{-1})}
\overline{g(t)}d\mu(t)=
\int_G\overline{f(t^{-1})}\overline{g(t)}d\mu(t^{-1})
$$
$$
=\int_G\overline{g(t^{-1})}\overline{f(t)}d\mu(t)=(S^*g,f),
$$
hence $(S^*g)(t)=\overline{g(t^{-1})}.$ Finally the modular operator
$\Delta$ defined by $\Delta=S^*S$ has the following form $(\Delta
f)(t)=\frac{d\mu(t)}{d\mu(t^{-1})}f(t)$. Indeed we have
$$
f(t) \stackrel{S}{\mapsto}
\frac{d\mu(t^{-1})}{d\mu(t)}\overline{f(t^{-1})}\stackrel{S^*}{\mapsto}
\frac{d\mu(t)}{d\mu(t^{-1})}f(t).
$$
Finally, since $J=S\Delta^{-1/2}$ (see \cite{{Con94}} p.462) we get
$$
f(t)
\stackrel{\Delta^{-1/2}}{\mapsto}\left(\frac{d\mu(t^{-1})}{d\mu(t)}\right)^{1/2}
f(t) \stackrel{S}{\mapsto}\frac{d\mu(t^{-1})}{d\mu(t)}
\left(\frac{d\mu(t)}{d\mu(t^{-1})}\right)^{1/2}\overline{f(t^{-1})}
$$
$$
=
\left(\frac{d\mu(t^{-1})}{d\mu(t)}\right)^{1/2}\overline{f(t^{-1})},
$$
\begin{equation}
\label{J,Delta-n}
(Jf)(t)=\left(\frac{d\mu(t^{-1})}{d\mu(t)}\right)^{1/2}\overline{f(t^{-1})},\,\,
\text{\,and\,}\,\, (\Delta f)(t)=\frac{d\mu(t)}{d\mu(t^{-1})}f(t).
\end{equation}
To prove that $JT^{R,\mu}_tJ=T^{L,\mu}_t$ we get
$$
f(t)\stackrel{J}{\mapsto}\left(\frac{d\mu(x^{-1})}{d\mu(x)}\right)^{1/2}
\overline{f(x^{-1})} \stackrel{T^{R,\mu}_t}{\mapsto}
\left(\frac{d\mu(xt)}{d\mu(x)}\right)^{1/2}
\left(\frac{d\mu((xt)^{-1})}{d\mu(xt)}\right)^{1/2}\overline{f((xt)^{-1})}
$$
$$
=\left( \frac{d\mu(t^{-1}x^{-1})}{d\mu(x)}
\right)^{1/2}\overline{f(t^{-1}x^{-1})}\stackrel{J}{\mapsto} \left(
\frac{d\mu(x^{-1})}{d\mu(x)} \right)^{1/2}\left(
\frac{d\mu(t^{-1}x)}{d\mu(x^{-1})} \right)^{1/2}f(t^{-1}x)
$$
$$
=\left(\frac{d\mu(t^{-1}x)}{d\mu(x)} \right)^{1/2} f(t^{-1}x)=
(T^{L,\mu}_tf)(x).
$$
\begin{rem} The representation $T^{R,\mu_b}$ is the inductive
limit of the representations $T^{R,\mu_b^m}$ of the group
$B(m,{\mathbb R})$ where the measure $\mu_b^m$ is the projection of
the measure $\mu_b$ onto subgroup   $B(m,{\mathbb R})$. Obviously
$\mu_b^m$ is equivalent with the Haar measure $h_m$ on $B(m,{\mathbb
R})$.
\end{rem}
\section{Structure of von Neumann algebras and the flow of weights
invariant of Connes and Takesaki}
Let us denote as before $M={\mathfrak A}^{R,b}=(T^{R,b}_s \mid s \in
B^{\mathbb Z})^{\prime\prime},\,$ ${\mathfrak A}^{L,b}=(T^{L,b}_t
\mid t \in B_0^{\mathbb Z})^{\prime\prime}$. We assume that
$$
E(b)<\infty, \text{\,\,\,
hence\quad}S^R_{kn}(b)<\infty,\,\,\,\text{and}\,\,\,
S^L_{kn}(b)<\infty.
$$
We also note that from \cite{Kos01m2}, Lemma 4 (see Remark
\ref{rem_ergo}) follows that the condition $E(b)<\infty$ implies the
ergodicity of the measure $\mu_b$ with respect to the right action.\\
From the previous section (see (\ref{J,Delta-n})) we conclude that
the modular operator $\Delta$ is defined as follows
\begin{equation}
\label{Delta-Z} (\Delta f)(x)=d\mu_b(x)/d\mu_b(x^{-1})f(x).
\end{equation}
In the next section we shall prove that $M$ (and hence $M^\prime$)
is a type III$_1$ factor, without assuming further conditions on the
measure. From our proof immediately follows that $M$ is a factor and
we do not use the conditions (\ref{S^(R,L)_(kn)(b)}).
In \cite{Kos08arIII.1} the case of the group
$B_0^{\mathbb N}$ was studied. There the author proved the type
III$_1$ property by directly showing that the fixed point algebra of
$M$ w.r.t. the modular group is trivial. The idea was to prove that
the one-parameter groups generated by multiplication operators by
the independent variables $x_{kn}$ are contained in the commutant of
the fixed point algebra. This, together with the fact that also all
translations $T^{R,b}_t, t\in B_0^{\mathbb N}$ are in the latter
commutant, implies the triviality of this fixed point algebra. In
the case of $B_0^{\mathbb Z}$, however, this method does not
directly give the answer, because of the presence of infinite sums
(see later). In this paper, we came up with a different approach, by
making use of the \emph{flow of weights} of a type III factor.
\index{flow of weights} \index{}

After the original {\it classification of type III factors} by
Connes in 1973 (\cite{Con73}), there came an equivalent
classification of type III factors using the flow of weights
invariant, which was concluded in the joint work of Connes and
Takesaki \cite{ConTa77}. The definition of the flow of weights
relies on the {\it duality theory for von Neumann algebras}, which
was discovered by Connes (\cite{Con73}) and Takesaki (\cite{Tak73}).
\begin{rem} (See \cite{Tak02}, chap.X, \S 2)
To a von Neumann algebra $M$ with an action $\sigma$ of a locally
compact abelian group $G$ one can associate another von Neumann
algebra $N$ and an action $\theta$ of the dual group $\hat G$. The
pair $(N,\theta)$ is called the \emph{dual} system.
\end{rem}
\index{theory!duality for von Neumann algebras}
For the case when $M$ is a type III factor and $\sigma$ is the {\it
modular automorphism group} of a {\it faithful semi-finite normal
weight}, this becomes an invariant, called \emph{the non-commutative
flow of weights}. The pair $(\mathcal C_N,\theta)$, where $\mathcal
C_N$ is the center of $N$, is called \emph{the flow of weights}. In
this section we review some basic facts about the structure of von
Neumann algebras and the flow of weights of a type III factor. For
more details we refer to e.g. \cite{Tak02}.
\index{group!automorphism!modular} \index{flow of
weights!non-commutative} \index{weight!faithful!semi-finite normal }
\begin{df}
{\it A ${\rm W}^*$-dynamical system} is a triple $(M,\alpha,G)$,
where $M$ is a von Neumann algebra, $G$ a locally compact group and
$\alpha$ an action of $G$ on $M$ by automorphisms, i.e. a strongly
continuous homomorphism from $G$ into $Aut(M)$.
\end{df}
\index{${\rm W}^*$-dynamical system} \index{crossed product}
A {\it crossed product} of $(M,\alpha,G)$ is defined as follows.
\begin{df}
Consider the following representations of $M$ and $G$ on
$L^2(G,\mathcal H)=L^2(G)\otimes\mathcal H$, where $\mathcal H$ is
the Hilbert space $M$ acts on:
$$
(\pi_\alpha(x)\xi)(s)=\alpha_s^{-1}(x)\xi(s),\ \ \ x\in M \ \ \ s\in
G,
$$
$$
(\lambda(t)\xi)(s)=\xi(t^{-1}s), \ \ \ s,t\in G \ \ \ \xi\in
L^2(G,\mathcal H).
$$
The representation $(\pi_\alpha,\lambda)$ is covariant, i.e.
$\pi_\alpha\circ\alpha_s(x)=\lambda(s)\pi_\alpha(x)\lambda(s)^*$.
Then
$$
M \rtimes _\alpha G:=(\pi_\alpha(M)\cup\lambda(G))^{\prime\prime}
$$
is called \emph{the crossed product} of $(M,\alpha,G)$.
\end{df}
When $G$ is abelian we can define an action of $\hat G$ (the dual of
$G$) on the crossed product by the following formulas.
$$
(\mu(p)\xi)(s)=\overline{\langle s,p\rangle}\xi(s),\ \ \ p\in \hat
G,
$$
$$
\hat\alpha_p(x)=\mu(p)x\mu(p)^*, \ \ \ x\in M \rtimes_\alpha G.
$$
Let us denote $U(s)=\lambda(s)$ and $V(p)=\mu(p)$, $s\in G, p\in
\hat G$. It
 follows that $U$ and $V$ satisfy the following relation:
\begin{equation}\label{eq_WH}
    U(s)V(p)U(s)^*V(p)^*=\langle s,p\rangle.
\end{equation}
\begin{df}
    In general, a pair of unitary representations $U$ of $G$ and $V$ of $\hat G$
     on the same Hilbert space $\mathcal H$ is said to be \emph{covariant} if
     the commutation relation (\ref{eq_WH}) is satisfied. The commutation
     relation (\ref{eq_WH}) is called the \emph{Weyl-Heisenberg} commutation
     relation.
\end{df}
\index{commutation relation!Weyl-Heisenberg}
 In what follows we will need the following result for
covariant representations in a von Neumann algebra.
\begin{prop}[\cite{Tak02}, Proposition 2.2]\label{prop_WH}
    The covariant representation $\{\lambda_G,\mu_G\}$ generates the factor
    $B(L^2(G))$ of all bounded operators on $L^2(G)$.
\end{prop}
\begin{pf}
 For each $f\in L^1(\hat G)$, we define
    \begin{equation*}
        V(f):=\int_{\hat G} f(p)V(p)dp.
    \end{equation*}
    Then $V$ is a *-representation of $L^1(\hat G)$, so that it can be extended
    to the enveloping $C^*$-algebra $C_0(G)$ (the algebra of continuous functions
    vanishing at infinity\footnote{A function $f$ on $G$ is said to vanish at
    infinity if given any $\epsilon>0$, there is a compact subset of $G$ such
    that $|f(x)|<\epsilon$ for $x$ outside this subset}). We shall denote the
    extended representation of $C_0(G)$ by $V$ again. In the case when $V=\mu_G$,
     we have that $\mu_G(f)$ is the multiplication by $f$ on $L^2(G)$
     ($f\in C_0(G)$). Hence the von Neumann algebra $A$ generated by
     $\{\mu_G(f); f\in C_0(G)\}$ is the multiplication algebra $L^\infty(G)$ on
      $L^2(G)$. So it is maximal abelian (i.e. $L^\infty(G)^\prime=L^\infty(G)$).
      Now, we have
    \begin{equation*}
        \lambda_G(s)\mu_G(f)\lambda_G(s)^*=\mu_G(\lambda_s f), \quad s\in G,
        \quad f\in L^\infty(G),
    \end{equation*}
    where $(\lambda_sf)(r)=f(r-s)$. Hence the operators  $A$ commuting
    with $\lambda_G(G)$ are only scalars (the Haar measure $dr$ is ergodic).
    Therefore,
\begin{equation*}
        \{\lambda_G(G), \mu_G(\hat G)\}^\prime=\mathbb C,
\end{equation*}
 so that $\{\lambda_G,\mu_G\}$ is irreducible.
\qed\end{pf}
The definition of the flow of weight relies on the following duality
theorem of Connes and Takesaki.
\begin{thm}[\cite{Con73,Tak73}]
For a W$^*$-dynamical system $(M,\alpha,G)$, where $G$ is abelian
the following isomorphism holds:
$$
(M \rtimes_\alpha G) \rtimes_{\hat\alpha}\hat G\cong
M\overline\otimes B(L^2(G)).
$$
\end{thm}
\begin{df}\label{df_dual}
    The representation $\mu$ of $\hat G$ on $L^2(G,\mathcal H)$ defined above
    is called \emph{the dual representation} to $\lambda$. The action
    $\hat\alpha$ of $\hat G$ on the crossed product $\hat M=M\rtimes_\alpha G$
     is called \emph{the dual action} and the resulting dynamical system
     $(\hat M,\hat\alpha,\hat G)$, we call \emph{the dual system}.
\end{df}
Later we will need a convenient description of the commutant of
$\hat M$. The following theorem gives us the desired answer:
\begin{thm}[\cite{Tak02}]\label{thm_Nprime}
    Consider a $W^*$-dynamical system $(M,G,\alpha)$, where $M$ acts on
    $\mathcal H$, $G$ a locally compact group and $\alpha$ is implemented by a
     unitary one parameter group $V(s)$ on $\mathcal H$, i.e.
     $\alpha_s(a)=V(s)aV(s)^*$, for $a\in M,s\in G$. Define
    \begin{equation*}
        (W\xi)(s)=V(s)^*\xi(s),\quad \xi\in L^2(G,\mathcal H).
    \end{equation*}
Then the following holds
$$M\rtimes_\alpha G=\left(WMW^*\cup \mathfrak A^\lambda_G\right)^{\prime\prime},$$
$$(M\rtimes_\alpha G)^\prime=\left(M^\prime\cup W\mathfrak A^\rho_GW^*\right)^{\prime\prime},$$
where $\mathfrak A^\lambda_G$ (resp. $\mathfrak A^\rho_G$) is the
left (resp. the right) von Neumann algebra of $G$.
\end{thm}
The following Theorem describes the so-called \emph{continuous
decomposition} of a von Neumann algebra. It is crucial for the
definition of flow of weights for type III factors. From now on we
set $G=\mathbb R$ and denote the $W^*$-dynamial systen
$(N,\theta,\mathbb R)$ by $(N,\theta)$.
\begin{thm}[\cite{Tak73,ConTa77}]\label{thm_structTypeIII}
 \begin{enumerate}
    \item   Let $(N,\theta)$ be a $W^*$-dynamical system such that
    \begin{itemize}
        \item $N$ admits a faithful semi-finite normal trace $\tau$;
        \item $\theta$ transforms in such a way that
        \begin{equation*}
            \tau\circ\theta_s=e^{-s}\tau, \quad s\in\mathbb R.
        \end{equation*}
    \end{itemize}
    Then the crossed product $M=N\rtimes_\theta \mathbb R$ is properly infinite
     and the center $\mathcal C_M$ is precisely the fixed point algebra
     $\mathcal C^\theta_N$ of the center of $N$ under the canonical embedding of
      $N$ into $M$ (the representation $\pi_\theta$). Furthermore, $M$ is of type
       III (i.e. all the factors in the decomposition of $M$ are of type III) if
        and only if the central dynamical system $(\mathcal C_N,\theta)$ does not
         contain an invariant subalgebra $\mathcal A$, such that the subsystem
         $(\mathcal A,\theta)$ is isomorphic to $L^\infty(\mathbb R)$ together
          with the translation action of $\mathbb R$. In the case that $M$ is of
           type III, $N$ is necessarily of type II$_\infty$
           (i.e. $\tau(I)=\infty$).
    \item If $M$ is a von Neumann algebra of type III, then there exists a unique,
     up to conjugacy, covariant system $(N,\theta)$ satisfying the conditions of
     (1).
    \end{enumerate}
\end{thm}
Moreover, the above theorem implies that $M$ is a factor if and only
if $\theta$ is {\it centrally ergodic}. Now we are ready to
introduce the flow of weights. It was discovered, in the context of
the duality theory for von Neumann algebras, by Connes in
\cite{Con73} and Takesaki in \cite{Tak73} and was studied in detail
in their joint work \cite{ConTa77}.
\index{$\theta$ centrally ergodic}
\begin{df}
The dynamical system $(N,\theta,\mathbb R)$, such that $M\cong N
\rtimes_\theta\mathbb R$ is called \emph{the non-commutative flow of
weights}, whereas $(\mathcal C_N,\theta,\mathbb R)$ is called
\emph{the flow of weights} associated to $(M,\sigma,\mathbb R)$. By
the above theorem it is an invariant for the algebraic type of $M$.
\end{df}
Recall that in \cite{Con73} Alain Connes classified type III factors
with the following invariant.
\begin{equation}
\label{S(M)-Z} S(M)=\bigcap_{\phi\in \mathcal W}Sp\Delta_{\phi},
\end{equation}
where $\mathcal W$ is the set of {\it faithful normal semi-finite
weights} on $M$. $S(M)\backslash \{0\}$ is a multiplicative subgroup
of $\mathbb R_+$ and subdivides type III factors into three classes
(\cite{Con73}):
\index{weight!faithful!normal semi-finite}
\begin{enumerate}
\item When $S(M)=[0,+\infty)$, $M$ is said to be of type III$_1$
\item When $S(M)=\{\lambda^n|n\in\mathbb Z\}\cup \{0\},$ where $0<\lambda<1$,
$M$ is said to be of type III$_\lambda$,
\item When $S(M)=\{0,1\},$ $M$ is said to be of type III$_0$.
\end{enumerate}
The following theorem states an equivalent description of the types
in terms
 of the flow of weights.
\begin{thm}[\cite{ConTa77}]\label{thm_flowofweights}
Let $M$ be a factor of type III.
\begin{enumerate}
\item M is of type III$_1$ if and only if the flow of weights is trivial, i.e.
$N$ is a factor.
\item $M$ is of type III$_0$ if and only if $N$ is not a factor and the flow of
 weights has no period.
\item $M$ is of type III$_\lambda$ if and only if $N$ is not a factor and $T>0$
is the period of the flow of weights,where $\lambda={\rm e}^{-T}$.
\end{enumerate}
\end{thm}
\section{Type III$_1$ factor}
For the von Neumann algebra $M=\mathfrak A^{R,b}$  we shall prove
that the corresponding flow of weights is trivial , i.e. $N$ is a
factor.
Using the theorem \ref{thm_flowofweights}, we conclude that $M$ is
then of type III$_1$.
[From Theorem \ref{thm_structTypeIII} it follows, first of all that
$M$ is a factor, since the center of $M$ is contained in the center
of its dual. Moreover, by the same theorem, it is of type III, since
there can not be a non trivial subsystem isomorphic to $L^2(\mathbb
R)$
 with the translation action on $\mathbb R$. Furthermore, by theorem
 \ref{thm_flowofweights}, $M$ is of type III$_1$.]\\
  Note that to prove the factor
 property we do not use the sufficient conditions from Corollary \ref{cor_factor}.
 Now we state the main theorem.
\begin{thm}\label{thm_main2}
Consider the von Neumann algebra ${\mathfrak A}^{R,b}$ generated by
the right regular representation $T^{R,b}$ of $B_0^{\mathbb Z}$.
Assume that $E(b)<\infty$. Let $\phi(a)=(1,a1)$ be the faithful
normal state on ${\mathfrak A}^{R,b}$, associated to the cyclic and
 separating
 vector $1$, and $\sigma$ the corresponding modular automorphism group. Then the
  dual algebra $N:={\mathfrak A}^{R,b}\rtimes_\sigma\mathbb R$ is a factor. Hence,
  ${\mathfrak A}^{R,b}$ is a
  type III$_1$ factor. The same holds for ${\mathfrak A}^{L,b}$.
\end{thm}
\begin{pf} On the space $L^2(\mathbb R,L^2(B^{\mathbb Z},\mu_b))=
L^2(B^{\mathbb Z}\times\mathbb R,\mu_b\otimes m)$ ($m$ is the
Lebesgue measure)   we define the operator $W$ as follows
\begin{equation}
\label{W=D^(-it)}
    (Wf)(x,t)=\Delta^{-it}(x)f(x,t)=\left(\frac{d\mu_b(x^{-1})}
    {d\mu_b(x)}\right)^{it}f(x,t).
\end{equation}
Denote $N:=\hat M=M\rtimes_\sigma\mathbb R$. By Theorem
\ref{thm_Nprime}, we have
$$
N=(W MW^*\cup {\mathfrak A}^\lambda_{\mathbb R}),\quad N^\prime=(
M^\prime\cup W{\mathfrak A}^\rho_{\mathbb R}W^*),
$$
hence
\begin{equation}\label{eq_CNprime}
    \mathcal C_N=N^\prime\cap N=(N\cup N^\prime)^\prime=
    (W MW^*\cup M^\prime\cup\lambda(\mathbb R)\cup W\rho(\mathbb
    R)W^*)^\prime.
\end{equation}
%
%
From (\ref{eq_CNprime}) we see that $\mathcal C_N^\prime$ contains
the following set of elements: 
\begin{equation}\label{eq_CNprim}
\left(WT^{R,b}_uW^*,T^{L,b}_u,\lambda(s),W\rho(s)W^*;u\in
B_0^{\mathbb Z},\,\,s\in{\mathbb R}\right).
\end{equation}
 We would like to prove the triviality of $\mathcal C_N$. For this  we show that
operators of multiplication by the independent variables $x_{kn}$
and $t$, in the space $L^2(B^{\mathbb Z}\times\mathbb R,\mu_b\otimes
m)$, are affiliated to $\mathcal C_N^\prime$ (see Lemma
\ref{lem_affil2}). In this case, since $T^{L,b}_u\in \mathcal
C_N^\prime, \forall u\in B_0^{\mathbb Z}$, the two lemmas below
would imply that $\mathcal C_N$ is trivial.
\begin{lem}[\cite{Kos08arIII.1}]
\label{lem_TgT}
    Let $g$ be a multiplication on $L^2(B^{\mathbb Z},\mu_b)$ by a measurable
    function $g$ on $B^{\mathbb Z}$, then
\begin{equation*}
    (T^{R,b}_tgT^{R,b}_{t^{-1}}f)(x)=g(xt)f(x), \quad\text{for all\,\,\,}
    t\in B_0^{\mathbb Z},f\in L^2(B^{\mathbb Z},\mu_b).
\end{equation*}
\end{lem}
%
%
\begin{lem}\label{prop_weylheis}
Let $M$ be a von Neumann algebra on $L^2(B^{\mathbb Z},\mu_b)\otimes
L^2(\mathbb R,m)= L^2(B^{\mathbb Z}\times\mathbb R,\mu_b\otimes m)$.
If $e^{ist}, e^{isx_{kn}}\in M^\prime, k<n$, $T^{L,b}_u\in
M^\prime,\forall u\in B_0^{\mathbb Z},s\in\mathbb R$, $\lambda(s)\in
M^\prime$ for all $s\in\mathbb R$ and the measure $\mu_b$ is
 ergodic, then $M={\mathbb C}I$.
\end{lem}
\begin{pf}
From  proposition \ref{prop_WH} follows that the result holds in the
one-dimensional case. The space $L^2(B^{\mathbb Z}\times\mathbb
R,\mu_b\otimes m)$ is isomorphic to $L^2(\mathbb
R^\infty,\mu_b\otimes m)= \otimes_{k<n\in\mathbb Z}L^2(\mathbb
R^1,\mu_{b_{kn}})\otimes L^2(\mathbb R,m)$.
 Since the variables $x_{kn}$ and $t$ are independent, the condition $e^{its},
 e^{ix_{kn}s}\in M^\prime$, for all $k<n\in \mathbb Z$ and $s\in\mathbb R$, means
 that $L^\infty(\mathbb R,m), L^\infty(\mathbb R,\mu_{b_{kn}})\subset M^\prime$
 for all $k<n$. This implies that the von Neumann algebra generated by
 $(L^\infty(\mathbb R,\mu_{b_{kn}}))_{k<n}\cup L^\infty(\mathbb R,m)$ is
 contained in $M^\prime$. The latter is isomorphic to
 $L^\infty(B^{\mathbb Z}\times\mathbb R,\mu_b\otimes m)$, which is maximally abelian.
 Hence, $M\subset L^\infty(B^{\mathbb Z}\times\mathbb R,\mu_b\otimes m)^\prime=
 L^\infty(B^{\mathbb Z}\times\mathbb R,\mu_b\otimes m)$. Moreover, since we assume that
$\lambda(s),T^{L,b}_u\in M^\prime$ for all $u\in B_0^{\mathbb Z},
s\in\mathbb R$, all functions in $M$ are $B_0^{\mathbb Z}$-left
invariant, by Lemma \ref{lem_TgT}, and translation invariant in the
last argument. By the ergodicity of the measure, they are constant
$\mu_b\otimes m$-a.e. 
\qed\end{pf} 
Thus, Theorem  \ref{thm_main2} is proved.\qed\end{pf}
\begin{lem}\label{lem_affil2}
Let $Q_{kn}$ and $Q_t$ be the multiplication operators
\begin{equation}\
\label{Q_kn,Q_t} (Q_{kn}f)(x,t):=x_{kn}f(x,t),\quad
(Q_tf)(x,t):=tf(x,t),
\end{equation}
in $L^2(B^{\mathbb Z}\times\mathbb R,\mu_b\otimes m)$. Then
$$
e^{iQ_{kn}s},e^{iQ_ts} \in \mathcal C_N^\prime,
$$
for all $s\in\mathbb R,\,\,\,k,n\in\mathbb Z,\,\,\,k<n$.
\end{lem}
{\bf Formal Computations}: The following method uses calculations
with unbounded generators of one parameter groups, using commutator
identities. In such a way we want to obtain the variables $x_{kn}$
and $t$. However, these computations are formal and we would have to
verify conditions on the domains of the operators in question, to
justify these identities rigorously.

Consider generators of the following one-parameter groups in
 $\mathcal C_N^\prime$:
\begin{equation}\label{gener-C'_N}
(WT^{R,b}_{1+sE_{pq}}W^*,T^{L,b}_{1+sE_{pq}},\lambda(t)\,|\,s,t\in\mathbb
R).
\end{equation}
Note that we left out one group, namely $W\rho(\mathbb R)W^*$.
However, the same
result can be obtained by replacing $\lambda(\mathbb R)$ by the latter group. \\
The generator of $(\lambda(-t))$ is $D_t:=\frac{d}{dt}$. From
\cite{KosZek01} follows that the generators of
$T^L_{pq}(s):=T^{L,b}_{1+sE_{pq}}$ are
\begin{equation}\label{A^L}
A^L_{pq}:=\frac{d}{ds}T^{L,b}_{1+sE_{pq}}\vert_{s=0}=\Sigma_{m=q+1}^\infty
x_{qm}D_{pm}+D_{pq},
\end{equation}
where $D_{pq}=\frac{\partial}{\partial x_{pq}}-b_{pq}x_{pq}$.
Finally we need to calculate the generators of
$V_{pq}(s):=WT^{R,b}_{1+sE_{pq}}W^*.$
We have
$$
(V_{pq}f)(x,t):=\frac{d}{ds}(V_{pq}(s)f)(x,t)\vert_{s=0}=
\frac{d}{ds}(WT^R_{pq}(s)\Delta W^*f)(x,t)\vert_{s=0}
$$
$$
=\Delta(x)^{-it}(\frac{d}{ds}\Delta(xs)^{-it}.{\bf
1})_{|s=0}f(x,t)+\frac{d}{ds}(T^R_{pq}(s)f)(x,t)\vert_{s=0}.
$$
The last term is nothing else than $A^R_{pq}$ defined by (see e.g.
\cite{KosZek00})
$$
A^{R}_{kn}=\sum_{r=-\infty}^{k-1}x_{kr}D_{rn}+D_{kn},\quad 1\leq
k<n.
$$
Set $B_{pq}=V_{pq}-\frac{d}{ds}T^R_{pq}(s)$. We have
$$
B_{pq}:=\frac{d}{ds}(\Delta(x)^{-it}\exp(sA^R_{pq})\Delta^{it}\exp(-sA^R_{pq}){\bf
1})(x)_{|s=0}
$$
$$
=\Delta^{-it}([A^R_{pq},\Delta^{it}]).
$$
The $nth$ term in the Taylor expansion of $\Delta^{it}$ is equal to
$\frac{(it)^n}{n!}\ln\Delta^n$. Since $[A^R_{pq},\ln\Delta]$ is a
function (see later), it commutes with $\ln\Delta$. Hence one has
the following formula:
$$
[A_{pq}^R,\frac{(it)^n}{n!}\ln\Delta^n]=n\ln\Delta^{n-1}[A^R_{pq},\ln\Delta].
$$
applying this to the Taylor expansion of $\Delta^{it}$ we obtain
$$
[A^R_{pq},\Delta^{it}]=it\Delta^{it}[A^R_{pq},\ln\Delta].
$$
In this manner we obtain
$$
B_{pq}=it[A^R_{pq},\ln\Delta(x)]
$$
and hence $V_{pq}=it[A^R_{pq},\ln\Delta(x)]+A^R_{pq}$. Thus the
generator in (\ref{gener-C'_N}) are as follows:
$$
V_{pq}:=it[A^R_{pq},\ln\Delta(x)]+A^R_{pq},\,\,\,A^L_{pq},\,\,\,D_t:=\frac{d}{dt}.
$$
{\bf Some useful formulas} (see  \cite{Kos92}).\\
 Let us denote by $X^{-1}$ the inverse matrix to the
upper triangular matrix $X=I+x=I+\sum_{k<n}x_{kn}E_{kn}\in
B^{\mathbb Z}$
$$
X^{-1}=(I+x)^{-1}:=I+\sum_{k<n}x_{kn}^{-1}E_{kn}\in B^{\mathbb Z}.
$$
We have by definition $X^{-1}X=XX^{-1}=I$, hence
\begin{equation}
\label{x_{kn}{-1}Z} \left(XX^{-1}\right)_{kn}=
\sum_{r=k}^{n}x_{kr}x_{rn}^{-1}=\delta_{kn}=\sum_{r=k}^{n}x_{kr}^{-1}x_{rn}
=\left(X^{-1}X\right)_{kn} ,\quad k\leq n,
\end{equation}
thus
$$
x_{kn}^{-1}+\sum_{r=k+1}^{n-1}x_{kr}x_{rn}^{-1}+x_{kn}=0=
x_{kn}+\sum_{r=k+1}^{n-1}x_{kr}x_{rn}^{-1}+x_{kn}^{-1},\quad k< n,
$$
and
\begin{equation}
\label{x{kn}(-1)1Z}
x_{kn}^{-1}=-x_{kn}-\sum_{r=k+1}^{n-1}x_{kr}x_{rn}^{-1}=
-x_{kn}-\sum_{r=k+1}^{n-1}x_{kr}^{-1}x_{rn}.
\end{equation}
We can write also
\begin{equation}
\label{x{kn}(-1)2Z} x_{kn}^{-1}=-\sum_{r=k+1}^{n}x_{kr}x_{rn}^{-1}=
-\sum_{r=k}^{n-1}x_{kr}^{-1}x_{rn}.
\end{equation}
There is also the explicit formula for $x_{kn}^{-1}$ (see
\cite{Kos88} formula (4.4)) \\$x_{kk+1}^{-1}=-x_{kk+1}$ and
\begin{equation}
\label{x{kn}(-1)Z}
x_{kn}^{-1}=-x_{kn}+\sum_{r=1}^{n-k-1}(-1)^{r+1}\sum_{k\leq
i_1<i_2<...<i_r\leq n }x_{ki_1}x_{i_1i_2}...x_{i_rn},\quad k<n-1.
\end{equation}
\begin{rem}
 Using (\ref{x{kn}(-1)Z}) we see that $x_{kn}^{-1}$ depends only on $x_{rs}$
 with $k\leq r<s\leq n$.
\end{rem}
 Using (\ref{x{kn}(-1)2Z})  we have
\begin{equation}
\label{(.+.)(.-.)Z}
x_{kn}^{-1}+x_{kn}=-\sum_{r=k+1}^{n-1}x_{kr}x_{rn}^{-1},\quad
x_{kn}^{-1}-x_{kn}=2x_{kn}-\sum_{r=k+1}^{n-1}x_{kr}x_{rn}^{-1}.
\end{equation}
Let us denote
$$
w_{kn}:=w_{kn}(x):=(x_{kn}+x_{kn}^{-1})(x_{kn}-x_{kn}^{-1}).
$$
 Using (\ref{mu-bZ}) and (\ref{Delta-Z}) we get
\begin{equation}
\label{Delta(x)Z} \Delta(x)=\frac{d\mu_b(x)}{d\mu_b(x^{-1})}
=\exp\left[\sum_{k,n\in{\mathbb
Z},\,k<n}b_{kn}\left((x_{kn}^{-1})^2-x_{kn}^2\right)
\right].
\end{equation}
$$
-\ln\Delta(x)=\sum_{k,n\in{\mathbb
Z},\,k<n}b_{kn}\left[x_{kn}^2-(x_{kn}^{-1})^2\right]=
\sum_{k,n\in{\mathbb
Z},\,k<n}b_{kn}(x_{kn}+x_{kn}^{-1})(x_{kn}-x_{kn}^{-1})
$$
$$
\sum_{k,n\in{\mathbb
Z},\,k<n}b_{kn}(x_{kn}+x_{kn}^{-1})[2x_{kn}-(x_{kn}+x_{kn}^{-1})]=
\sum_{k,n\in{\mathbb Z},\,k<n}b_{kn}w_{kn}(x).
$$
To study the action of the operators
$$
A^{R}_{kn}:=\frac{d}{ds}T^{R,b}_{1+sE_{pq}}\vert_{s=0}=\sum_{r=-\infty}^{k-1}x_{rk}D_{rn}+D_{kn}
$$
on the function $\ln\Delta(x)$ we need to know the action of
$D_{pq}$ on $x_{kn}^{-1}$. 
\begin{lem}[ \cite{Kos08arIII.1}]
\label{l.D,x^{-1}Z} We have
\begin{equation}
\label{[D,x{-1}]Z} [D_{pq},x_{kn}^{-1}]= \left\{\begin{array}{cc}
-x_{kp}^{-1}x_{qn}^{-1},&\text{\,if\,}\,k\leq p<q\leq n,\\
0,&\text{\,otherwise\,}.
\end{array}\right.
\end{equation}
\end{lem}
For proof see \cite{Kos08arIII.1}.\\
Using (\ref{[D,x{-1}]Z}) we get
\begin{equation}
\label{[D,x+x{-1}]Z} [D_{pq},x_{kn}+x_{kn}^{-1}]=
\left\{\begin{array}{cc}
-x_{kp}^{-1}x_{qn}^{-1}-\delta_{kp}x_{qn}^{-1}-\delta_{qn}x_{kp}^{-1},&
\text{\,if\,}\,k\leq p<q\leq n,\\
0,&\text{\,otherwise\,}.
\end{array}\right.
\end{equation}
Using (\ref{[D,x+x{-1}]Z}) we have
$[D_{pq},(x_{kn}+x_{kn}^{-1})(x_{kn}-x_{kn}^{-1})]=$
\begin{equation}
\label{[D,w]Z} \left\{\begin{array}{cc}
2\,x_{kp}^{-1}x_{qn}^{-1}x_{kn}^{-1}+2\,\delta_{kp}x_{qn}^{-1}x_{kn}^{-1}+2\,
\delta_{qn}x_{kp}^{-1}x_{kn}^{-1},&\text{\,if\,}\,k\leq p<q\leq
n,\,(p,q)\not=
(k,n),\\
2(x_{kn}+x_{kn}^{-1}),&\text{\,if\,}\,(p,q)=(k,n),\\
0,&\text{\,otherwise\,}.
\end{array}\right.
\end{equation}
Indeed, if $k\leq p<q\leq n,\,(p,q)\not=(k,n)$ we have
$$
[D_{pq},(x_{kn}+x_{kn}^{-1})(x_{kn}-x_{kn}^{-1})]=
[D_{pq},(x_{kn}+x_{kn}^{-1})(2x_{kn}-(x_{kn}+x_{kn}^{-1}))]
$$
$$
=[D_{pq},(x_{kn}+x_{kn}^{-1})](2x_{kn}-(x_{kn}+x_{kn}^{-1}))-
(x_{kn}+x_{kn}^{-1})[D_{pq},(x_{kn}+x_{kn}^{-1})]=
$$
$$
-2x_{kn}^{-1}[D_{pq},(x_{kn}+x_{kn}^{-1})]
\stackrel{(\ref{[D,x+x{-1}]Z} )}{=}
2x_{kp}^{-1}x_{qn}^{-1}x_{kn}^{-1}+2\delta_{kp}x_{qn}^{-1}x_{kn}^{-1}+
2\delta_{qn}x_{kp}^{-1}x_{kn}^{-1}.
$$
\begin{lem}[ \cite{Kos08arIII.1}]  We have
\begin{equation}
 \label{[A^R,w]Z}
[A^R_{pq} ,\,w_{kn}]= \left\{\begin{array}{ll}
2x_{kp}x_{kq},           &\text{\,if\,}\,\,\,n=q,\,\,k\leq p-1,\\
2x_{pn}^{-1}x_{qn}^{-1}, &\text{\,if\,}\,\,\,k=p,\,\, n\geq q+1,\\
2(x_{pq}+x_{pq}^{-1}),   &\text{\,if\,}\,\,\,\,\,k=p,q=n,\\
0,                       &\text{\,if\,}\,\,\,k\not=p\text{\,\,or\,\,}n\not=q,\\
\end{array}\right.
\end{equation}
hence
\begin{equation}
\label{[A^R,lnDe]Z}
 -[A^R_{pq}
,\,\ln\Delta(x)]=2\sum_{r=-\infty}^{p-1}b_{rq}x_{rp}x_{rq} +
2\sum_{n=q+1}^\infty
b_{pn}x_{pn}^{-1}x_{qn}^{-1}+2(x_{pq}+x_{pq}^{-1}).
\end{equation}
\end{lem}
Next, we consider the action of $A^L_{ij}$ on (\ref{[A^R,lnDe]Z}),
where $i<p<j<q$.
\begin{lem} [ \cite{Kos08arIII.1}] One has
$$
[A^L_{ij},[A^R_{pq},\ln\Delta(x)]]=-2b_{iq}x_{jq}x_{ip},\ \
\text{for } i<p<j<q
$$
and hence
\begin{equation}\label{x_tq}
[A^L_{ip},[A^L_{ij},[A^R_{pq},\ln\Delta(x)]]]=2b_{iq}x_{jq},
\end{equation}
which immediately gives us the variables $x_{jq}$, $j,q\in\mathbb
Z,\ q-j\geq 1$.
\end{lem}
\begin{pf}
Recall (see (\ref{A^L})) that  $A^L_{ij}=\Sigma_{m=j+1}^\infty
x_{jm}D_{im}+D_{ij}$. Then
\begin{equation}\label{[A^L,x_kn{-1}]}
\begin{array}{ll}
[A^L_{ij},x_{kn}^{-1}]=\Sigma_{m=j+1}^\infty
x_{jm}[D_{im},x_{kn}^{-1}]+[D_{ij},
x_{kn}^{-1}]\\
\stackrel{(\ref{[D,x{-1}]Z})}{=}-\sum_{r=j+1}^{n}x_{jr}(x_{ki}^{-1}x_{rn}^{-1}+
\delta_{ki}x_{rn}^{-1}+\delta_{rn}x_{ki}^{-1}+\delta_{ki}\delta_{rn})\\
=-\delta_{jn}(x_{ki}^{-1}+\delta_{ki}).
\end{array}
\end{equation}
Moreover we also need the formula
\begin{equation}\label{[A^L,x_kn]}
\begin{array}{ll}
[A^L_{ij},x_{kn}]=\Sigma_{m=i+1}^\infty x_{jm}[D_{im},x_{kn}]+[D_{ij},x_{kn}]\\
=-\delta_{ki}(x_{jn}+\delta_{jn}).
\end{array}
\end{equation}
By our choice of $i$ and $j$, $[A^L_{ij},[A^R_{pq},\ln\Delta(x)]]$
will only depend on $x_{kn}$, $n<p$. Hence
$$
[A^L_{ij},[A^R_{pq},\ln\Delta(x)]]=-2\sum_{r=-\infty}^{p-1}b_{rq}x_{rp}[A^L_{ij},
x_{rq}]=-2b_{iq}x_{jq}x_{ip}.
$$
The last formula of the Lemma follows trivially. \qed\end{pf}
By applying $[A^L_{jq},.]$ to (\ref{x_tq}), we obtain a constant
$2b_{iq}$. The direct computation, based on the above formulas gives
us the operator of multiplication by $x_{kn}$ and $t$:
$$
[A^L_{jq},[ A^L_{ip},[A^L_{ij},C^R_{pq}]]]=2itb_{iq},
$$
$$
[D_t,[A^L_{ip},[A^L_{ij},C^R_{pq}]]]=2ib_{iq}x_{jq},
$$
since $A^R$ commutes with $A^L$,  $t$ is of course the variable in
$L^2(\mathbb R)$.
\begin{rem} The previous manipulations with the {\it unbounded operators},
were formal. Nevertheless they indicate us the form of the
expressions in terms of the {\it unitary one-parameter groups}
generated by $A^L_{kn}$ and $A^R_{kn}$, we should take, to obtain
the desired answer. We should {\it replace the commutator}
$[x,y],\,\,x,y\in L(G)$ in the {\it Lie algebra} $L(G)$ by  {\it the
group commutator} $\{a,b\}:=aba^{-1}b^{-1}$ where $a,b\in G$.
\end{rem}
Again, denote
$$
T^{L,b}_{pq}(s)=T^{L,b}_{I+sE_{pq}},\quad
T^{R,b}_{pq}(s)=T^{R,b}_{I+sE_{pq}},
$$
$$
V_{pq}(s)=WT^{R,b}_{I+sE_{pq}}W^*,\quad
Wf(x,t)=\Delta^{-it}(x)f(x,t).
$$
\begin{lem}\label{l.U(tau,s)} Denote
\begin{equation}
U(\tau,s)=\{T^{L,b}_\tau,V^{R,b}_s\},\,\,\tau,s\in B_0^{\mathbb Z},
\end{equation}
then we have
\begin{equation}
\label{U(tau,s)=} U(\tau,s)=
\Delta^{-it}(\tau^{-1}x)\Delta^{it}(\tau^{-1}xs)
\Delta^{-it}(xs)\Delta^{it}(x).
\end{equation}
\end{lem}
\begin{pf} Since
$$
U(\tau,s)=\{T^{L,b}_\tau,V^{R,b}_s\}=T^{L,b}_\tau
V^{R,b}_sT^{L,b}_{\tau^{-1}} V^{R,b}_{s^{-1}}= T^{L,b}_\tau
WT^{R,b}_{s}W^*T^{L,b}_{\tau^{-1}}WT^{R,b}_{s^{-1}}W^*,
$$
%
%
\index{ $\{a,b\}$-group commutator}
%
%
we have
$$
f(x,t)\stackrel{W^*\,\,}{\rightarrow}\Delta^{it}(x)f(x,t)
\stackrel{WT^{R,b}_{s^{-1}}\,\,}{\rightarrow}
\left(\frac{d\mu(xs^{-1})}{d\mu(x)}\right)^{1/2}\Delta^{it}(xs^{-1})f(xs^{-1},t)
\stackrel{W^*T^{L,b}_{\tau^{-1}}\,\,\,\,}{\rightarrow}
$$
$$
\Delta^{it}(x) \left(\frac{d\mu(\tau
x)}{d\mu(x)}\right)^{1/2}\Delta^{-it}(\tau x) \left(\frac{d\mu(\tau
xs^{-1})}{d\mu(\tau x)}\right)^{1/2}\Delta^{it}(\tau xs^{-1})f(\tau
xs^{-1},t)
$$
$$
=\Delta^{it}(x) \Delta^{-it}(\tau x) \left(\frac{d\mu(\tau
xs^{-1})}{d\mu(x)}\right)^{1/2}\Delta^{it}(\tau xs^{-1})f(\tau
xs^{-1},t)
$$
$$
\stackrel{WT^{R,b}_{s}\,\,}{\rightarrow}\Delta^{-it}(x)
\left(\frac{d\mu(\tau x)}{d\mu(x)}\right)^{1/2}\Delta^{-it}(xs)
\Delta^{-it}(\tau xs)\Delta^{it}(\tau x)f(\tau
x,t)\stackrel{T^{L,b}_{s\tau}\,\,}{\rightarrow}
$$
$$
\left(\frac{d\mu(\tau^{-1}x)}{d\mu(x)}\right)^{1/2}
\Delta^{-it}(\tau^{-1}
x)\left(\frac{d\mu(x)}{d\mu(\tau^{-1}x)}\right)^{1/2}
\Delta^{it}(\tau^{-1} xs)\Delta^{-it}(xs)\Delta^{-it}(x)f(x,t)
$$
$$
=\Delta^{-it}(\tau^{-1}x)\Delta^{it}(\tau^{-1}xs)
\Delta^{-it}(xs)\Delta^{it}(x)f(x,t).
$$
\qed\end{pf}
Consider the following one-parameter groups in $B_0^{\mathbb Z}$:
\begin{equation}\label{eq_Ekn(t)}
    E_{pq}(s):=\left\{1+sE_{pq}; s\in\mathbb R\right\},\,\, p,q\in{\mathbb
    Z},\,\,p<q.
\end{equation}
We calculate $U(E_{rm+1}(t),E_{mm+1}(s))$ for $t=-1$.
\begin{lem}\label{l.lem_U2}
Let $U_{rm}(s)\in \calC_N^\prime$ be the operators defined by
$$
U_{rm}(s):=U(I-E_{rm+1},I+sE_{mm+1}),
$$
where $s\in\mathbb R$. Then for $f\in L^2(B^{\mathbb Z}\times\mathbb
R,\mu_b\otimes m)$ holds
\begin{equation}\label{eq_Urm}
\left(U_{rm}(s)f\right)(x,t)=\exp\left(-2ib_{rm+1}stx_{rm}\right)f(x,t),\,\,\,\forall
t,s\in \mathbb R.
\end{equation}
Then  $\{U_{rm}(s),\lambda(1)\}$ and $\{U_{rm}(s),T^{L,b}_{rm}(1)\}$
are the one parameter groups generated by the multiplications by
$x_{rm}$ and $t$.
\end{lem}
\begin{pf}
Fix $s\in\mathbb R$ and define:
$$
X^\prime:=XE_{mm+1}(s),\quad Y:=E_{rm+1}(1)X,\quad
Y^\prime:=YE_{mm+1}(s),$$ where $X=1+x, Y=1+y\in B^{\mathbb Z}$.
First we note that $E_{rm+1}(-1)Y=X$ and
$E_{rm+1}(-1)Y^\prime=X^\prime$, since $E_{rm+1}(s)$ are
one-parameter groups.
By Lemma \ref{l.U(tau,s)} we get
$$
U_{rm}(s)=\frac{\Delta^{it}((I+E_{rm+1})x(I+sE_{mm+1}))
\Delta^{it}(x)}{\Delta^{it}((I+E_{rm+1})x)\Delta^{it}(x(I+sE_{mm+1}))},\text{\,\,\,\,or}
$$
$$
U_{rm}(s)=\Delta^{-it}(y)\Delta^{it}(y^\prime)
\Delta^{-it}(x^\prime)\Delta^{it}(x).
$$
To obtain (\ref{eq_Urm}) we proceed in two steps. First of all we
compute $\Delta(x^\prime)^{-it}\Delta(x)^{it}$. Secondly, we replace
$x$ by $y$ and $t$ by $-t$ in the latter expression, to obtain
$\Delta(y^\prime)^{it}\Delta(y)^{-it}$. Finally, we combine the two
expressions above to obtain the desired result.
 Recall that
$$
-\ln\Delta(x)=\sum_{k,n\in\mathbb Z,k<n}b_{kn}w_{kn}(x).
$$
We show (\ref{eq_Urm}) in two steps.\\
{\bf Step 1:} The first step  is to compute
$\Delta(x^\prime)^{-it}\Delta(x)^{it}$.
\begin{rem}
The computations below are analogous to those in
\cite{Kos08arIII.1}, which were carried out to obtain the variables
$x_{kn}$, in order to prove the triviality of the fixed point
algebra of ${\mathfrak A}^{L,b}$ w.r.t. the modular group. However,
in the case of $B_0^{\mathbb Z}$ the fixed point algebra method does
not work, since all sums are over an infinite number of indices.
 The \emph{flow of weights} allows us to \emph{overcome the problems}.
\end{rem}
First of all, we would like to know the right action of
$E_{mm+1}(s)$ on $X\in B^{\mathbb Z}$, where $X=1+x$ and on
$X^{-1}$. For $k<n$ we have: $x^\prime_{kn}=$
\begin{equation}\label{eq_xprime}
\left(XE_{mm+1}(s)\right)_{kn}=\sum_{i=k}^{n}X_{ki}(\delta_{in}+s\delta_{mi}
\delta_{m+1n})=x_{kn}+\delta_{m+1n}(sx_{km}+s\delta_{km}).
\end{equation}
Note that only the $m+1$st column of $x$ is affected by this
transformation.
\begin{equation*}
    XE_{mm+1}(s)=\left(
    \begin{array}{cccccccc}
        \ddots&&&&&\\
        &1&x_{12}&x_{13}&x_{14}&\cdots&x_{1m+1}+sx_{1m}&\cdots\\
        &0&1&x_{23}&x_{24}&\cdots& x_{2m+1}+sx_{2m}&\cdots\\
        &\vdots&\vdots&\vdots&\vdots&\vdots&\vdots&\vdots\\
        &0&0&\cdots&\cdots&1&x_{mm+1}+s&\cdots\\
        &&&&&&                      & \ddots
    \end{array}
    \right).
\end{equation*}
From (\ref{eq_xprime}) we now can deduce the inverse of $x^\prime$.
In order to do this we note the following:
$$x^{\prime-1}_{kn}=(XE_{mm+1}(s))_{kn}^{-1}=(E_{mm+1}(-s)X^{-1})_{kn},$$
where $k<n$. Hence we get: $x^{\prime-1}_{kn}=\sum_{j=k}^n
(E_{mm+1}(-s))_{kj}X_{jn}^{-1}=$
$$
\sum_{j=k}^n (\delta_{kj}-s\delta_{mk}
\delta_{jm+1})(x_{jn}^{-1}+\delta_{jn})=x_{kn}^{-1}-s\delta_{km}(x_{m+1n}^{-1}+
\delta_{m+1n}),
$$
where of course $x^\prime$ depends on $s\in\mathbb R$. Since only
the rows with number $m$ of $x^{-1}$ and the columns with number
$m+1$ of $x$ are affected by the transformation $x\to x^\prime$, it
follows that $\ln\Delta(x^\prime)$ can be written as the sum of
three different terms:
$$
-\ln\Delta(x^\prime)=\sum_{k<n,k\neq m,n\neq
m+1}b_{kn}w_{kn}(x)+\sum_{k<m}b_{km+1}w_{km+1}(x^\prime)+\sum_{n>m+1}b_{mn}
w_{mn}(x^\prime).
$$
Note that the term with $k=m,n=m+1$ vanishes, because
$w_{mm+1}(x)=0$. First we consider the second term:
$\sum_{k<m}b_{km+1}w_{km+1}(x^\prime)=$
\begin{equation}\label{eq_term2}
\begin{array}{l}
\sum_{k=-\infty}^{m-1}b_{km+1}(x^\prime_{km+1}-x^{\prime-1}_{km+1})
(x^\prime_{km+1}+x^{\prime-1}_{km+1})=\\  \quad
\sum_{k=-\infty}^{m-1}b_{km+1}(x_{km+1}+sx_{km}-x^{-1}_{km+1})
(x_{km+1}+sx_{km}+x^{-1}_{km+1})=\\  \quad
\sum_{k=-\infty}^{m-1}b_{km+1}w_{km+1}(x)+
\sum_{k<m}b_{km+1}\left(2sx_{km}x_{km+1}+s^2(x_{km})^2\right).
\end{array}
\end{equation}
 The third term is as follows:
$\sum_{n>m+1}b_{mn}w_{mn}(x^\prime)=$
\begin{equation}\label{eq_term3}
\begin{array}{l}
\sum_{n=m+2}^{\infty}b_{mn}(x^\prime_{mn}-x^{\prime-1}_{mn})(x^\prime_{mn}+x^{\prime-1}_{mn})=\\
\quad \sum_{n=m+2}^{\infty}b_{mn}(x_{mn}-x^{-1}_{mn}+sx_{m+1n}^{-1})
(x_{mn}+x_{mn}^{-1}-sx_{m+1n}^{-1})=\\
 \quad
\sum_{n=m+2}^{\infty}b_{mn}w_{mn}(x)+\sum_{n>m+1}b_{mn}\left(2sx_{mn}^{-1}
x_{m+1n}^{-1}-s^2(x_{m+1n}^{-1})^2\right).
\end{array}
\end{equation}
After adding up the terms we get:
$\ln\Delta(x)-\ln\Delta(x^\prime)=$
\begin{equation*}
    \sum_{k-\infty}^{m-1}b_{km+1}\left(2sx_{km}x_{km+1}+s^2(x_{km})^2\right)+
    \sum_{n=m+2}^{\infty}b_{mn}
    \left(2sx_{mn}^{-1}x_{m+1n}^{-1}-s^2(x_{m+1n}^{-1})^2\right).
\end{equation*}
Hence we get:
\begin{equation}\label{eq_U22}
\begin{array}{l}
\quad
\Delta^{it}(x)\Delta^{-it}(x^\prime)=\exp\left(it\sum_{k=-\infty}^{m-1}b_{km+1}(2sx_{km}x_{km+1}+s^2(x_{km})^2
\right)\\  \quad
+it\sum_{n=m+2}^{\infty}b_{mn}\left(2sx_{mn}^{-1}x_{m+1n}^{-1}-s^2
(x_{m+1n}^{-1})^2\right).
\end{array}
\end{equation}
{\bf Step 2:} The next step is to compute
$\Delta^{-it}(y)\Delta^{it}(y^\prime)$. To obtain the latter we have
to perform the following operations in formula (\ref{eq_U22})
\begin{equation*}
    x\to E_{rm+1}(1)x,\ \ \ t\to -t.
\end{equation*}
First we have to compute: $y_{kn}= \left(E_{rm+1}(1)X\right)_{kn}$
\begin{equation}\label{eq_primex}
=\sum_{i=k+1}^{n-1}\left(\delta_{ki}+\delta_{kr}
\delta_{im+1}\right)(x_{in}+\delta_{in})=x_{kn}+\delta_{kr}\left(x_{m+1n}+\delta_{m+1n}
\right).
\end{equation}
In order to calculate $y^{-1}$ we note that
$\left(E_{rm+1}(s)X\right)^{-1} =X^{-1}E_{rm+1}(-s)$ ($E_{rm+1}(s)$
are one-parameter groups). Thus
\begin{equation}\label{eq_primexinv}
y^{-1}_{kn}=x_{kn}^{-1}-\delta_{m+1n}(x^{-1}_{km}+\delta_{km}).
\end{equation}
According to equation (\ref{eq_primex}) only the row with number $r$
of $x$ is affected by the left $E_{rm+1}(1)$-action. Similarly, by
(\ref{eq_primexinv}), only the column with number $m+1$ of $y^{-1}$
is affected. Moreover, $y_{rm}=x_{rm}$, since $x_{m+1m}=0$. Hence,
in the calculation below we have:
\begin{equation*}
  y_{km}y_{km+1}=\left\{\begin{array}{ll}
    x_{km}x_{km+1},&\text{ if } k\neq r\\
    x_{rm}x_{rm+1}+x_{rm},&\text{ if } k=r,
    \end{array}\right.
\end{equation*}
and other terms are unaffected by the transformation $x\to y$.
We obtain
\begin{align*}
\Delta^{-it}(y)\Delta^{it}(y^\prime)=\\
\exp\{-it\sum_{k=-\infty}^{m-1}b_{km+1}(2s y_{km}
y_{km+1}+s^2(y_{km})^2)-\\
it\sum_{n=m+2}^{\infty}b_{mn}(2s y_{mn}^{-1} y_{m+1n}^{-1}-s^2(y_{m+1n}^{-1})^2\}=\\
\exp\{-it\sum_{k=-\infty}^{m-1}b_{km+1}(2sx_{km}x_{km+1}+s^2(x_{km})^2)-\\
it\sum_{n=m+2}^{\infty}b_{mn}(2sx_{mn}^{-1}x_{m+1n}^{-1}-s^2
(x_{m+1n}^{-1})^2) -2isb_{rm+1}tx_{rm}\}.
\end{align*}
Therefore
$$
\left(U_{rm}(s)f\right)(x,t)=e^{-2ib_{rm+1}stx_{rm}}f(x,t).
$$
Since $T^{L,b}_{rm}(s),\,\lambda(s)\in \calC_N^\prime$ for all $s\in
\mathbb R$, we conclude that $\{U_{rm}(s),\lambda(1)\},$ $
\{U_{rm}(s),\,T^{L,b}_{rm}(1)\}\in  \calC_N^\prime$. The explicit
calculation gives us:
\begin{eqnarray*}
\left(\{U_{rm}(s),\lambda(1)\}f\right)(x,t)&=&e^{isb_{rm+1}
x_{rm}}f(x,t),\\ \\
\left(\{U_{rm}(s),T^{L,b}_{rm}(1)\}f\right)(x,t)&=&
e^{isb_{rm+1}t}f(x,t).
\end{eqnarray*}\qed\end{pf}
\begin{pf} Now we finish the proof of Lemma \ref{lem_affil2}.
 From the equations above we see that the unitary
one-parameters groups $\exp(isQ_{rm})$ and $\exp(isQ_t),\,\,s\in
{\mathbb R}$, generated by the self-adjoint  operators $Q_{rm}$ and
$Q_t$, defined by (\ref{Q_kn,Q_t}), are contained in
$\calC_N^\prime$, which proves Lemma \ref{lem_affil2}. \qed\end{pf}
\section{Uniqueness of the constructed factor}
\begin{thm}\label{thm_unique}
The von Neumann algebras $\mathfrak A^{R,b}$ and $\mathfrak A^{L,b}$
are hyperfinite type III$_1$ factors and hence isomorphic to the
factor $R_\infty$
 of Araki and Woods.
\end{thm}
\begin{pf}
Let G be a solvable separable locally compact group or a connected
locally compact group. Then any representation $\pi$ of $G$ in a
Hilbert space generates a hyperfinite von Neumann algebra (\cite{Con76}).\\
The group $B_0^\mathbb Z$ is the inductive limit of groups of finite
dimensional upper-triangular matrices (with units on the diagonal),
which are of course solvable, connected locally compact groups.
Hence their group algebras are hyperfinite (and by a theorem of
Dixmier (\cite{Dix59}) even type I algebras). Thus the von Neumann
algebra $\mathfrak A^{R,b}$ is the inductive limit of hyperfinite
von Neumann algebras and hence itself hyperfinite. From the theorem
of Haagerup (\cite{Haa87}) follows that ${\mathfrak A}^{R,b}$ and
${\mathfrak A}^{L,b}$ are all isomorphic to the Araki-Woods factor
$R_\infty$ (\cite{ArW76}). \qed\end{pf}

{\it Acknowledgements.} {The second author would like to thank the
Max-Planck-Institute of Mathematics and the Institute of Applied
Mathematics, University of Bonn for the hospitality. The partial
financial support by the DFG project 436 UKR 113/87 are gratefully
acknowledged.

The first author was supported by the Max-Planck Doctoral
scholarship. He would like to thank the Max-Planck-Institute for
Mathematics for its hospitality and support. He is grateful to M.
Marcolli, D. Markiewicz, U. Haagerup, R. Longo, G. Elliottand D.
Guido for useful discussions. Parts of this project was done by the
first author during his visits to the Fields Institute as a
participant of the thematic program on Operator Algebras and to
Florida State University. He would like to thank the Fields
Institute and FSU for their hospitality and financial support.

Both authors are grateful to A.~Connes for  useful remarks
concerning uniqueness of the constructed factor.}

\end{document}